\documentclass[11pt,twoside]{amsart}
\usepackage[francais,english]{babel} 
\usepackage[T1]{fontenc}
\usepackage[utf8]{inputenc}
\usepackage{lmodern}
\usepackage{amsmath,amscd}
\usepackage{amsfonts}
\usepackage{amssymb}
\usepackage{amsthm}
\usepackage{newlfont}
\usepackage{graphicx}
\usepackage{float}
\usepackage{hyperref}

\theoremstyle{plain}
\newtheorem{thm}{Theorem}[section]
\newtheorem{cor}[thm]{Corollary}

\newtheorem{prop}[thm]{Proposition}
\theoremstyle{definition}
\newtheorem{defn}[thm]{Definition}
\theoremstyle{remark}
\newtheorem{rem}[thm]{Remark}
\theoremstyle{Question}
\newtheorem{quest}[thm]{Question}
\theoremstyle{Example}
\newtheorem{ex}[thm]{Example}
\theoremstyle{Problem}

\begin{document}

\title{Positive loops of loose Legendrian embeddings and applications}

\author{Guogang LIU}
\email{guogangliu2016@gmail.com}  
\address{Laboratoire de Math\'ematiques Jean Leray  \\
Universit\'e de Nantes\\
2 Chemin de la Houssiniere, 44322 Nantes\\
France}


\date{}

\maketitle

\begin{abstract}
  In this paper,  we prove that there exist contractible
 positive loops of  Legendrian  embeddings
 based at any loose Legendrian submanifold. As an application, we define a partial order on $\widetilde{Cont}_0(M,\xi)$, called strong orderability, and prove that  overtwisted contact manifolds are
 not strongly orderable.
\end{abstract}
\tableofcontents

\section*{Introduction}

In this paper, we focus on the study of positive contact and Legendrian isotopies in a 
 co-oriented contact manifold $(M,\xi)$.

A contact manifold $(M^{2n+1},\xi)$ is a 
$2n+1$ dimensional smooth manifold $M$ with
a non-integrable hyperplane field $\xi$ which
is called a contact structure.
When $\xi$ is co-oriented, it is given by the 
kernel of a  \textit{contact} $1$-form $\alpha$. For example,
in $\mathbb{R}^4$ with the usual coordinates
$(x_1,x_1,y_1,y_2)$, the sphere $\mathbb{S}^3$
carries a contact form $\alpha_{std}=(y_1dx_1-x_1dy_1 +y_2dx_2-x_2dy_2 )|_{\mathbb{S}^3}$. We denote $\xi_{std}$ the 
contact structure defined by $\alpha_{std}$.
It induces a contact structure on 
 the quotient $\mathbb{R}P^3$  which is also denoted
 by $\xi_{std}$.
 
 One class of  submanifolds of $(M^{2n+1},\xi)$ with an interesting behavior is that of Legendrian submanifolds. A $n$-dimensional
 submanifold $L\subset M^{2n+1}$
  is  called a Legendrian submanifold if $\alpha|_L=0$.
A contactomorphism of $(M,\xi)$ is a diffeomorphism
which preserves $\xi$ and a contact isotopy 
$(\varphi_t)_{t\in [0,1]}$ is a path of contactomorphisms with $\varphi_0=id$. We say a contact isotopy $(\varphi_t)_{t\in [0,1]}$
is positive if $\alpha(\partial_t\phi_t)>0$.
That is to say, the infinitesimal generator of the isotopy is 
positively transverse to $\xi$ everywhere. An isotopy 
$(\varphi_t)_{t\in [0,1]}$ based at a Legendrian submanifold $L$
is said to be a Legendrian isotopy if $\varphi_t(L)$
is a Legendrian submanifold for all $t$. Similarly,
We say $\varphi_t$ is positive if $\alpha(\partial_t\varphi_t)>0$. This notion of positivity does only depend on the image $L_t =\varphi_t (L)$ of the isotopy. For us, a Legendrian isotopy will be such a family of unparametrized Legendrian submanifolds.

With the concept of positive contact isotopy, 
Eliashberg and Polterovich defined a partial order on
the universal cover $\widetilde{Cont}_0(M,\xi)$ of the identity component of the contactomorphisms group of $(M,\xi )$.
A class of contact isotopy $[(\psi_t)_{t\in [0,1]}]$
is greater than another class $[(\varphi_t)_{t\in [0,1]}]$
if there exists a positive contact isotopy from
$\varphi_1$ to $\psi_1$ which is homotopic to the concatenation
of the opposite of $(\varphi_t)_{t\in [0,1]}$ 
and $(\psi_t)_{t\in [0,1]}$.

\begin{prop}\cite{EP99}\label{EP99}
If $(M, \xi)$ is a contact manifold, the following conditions
are equivalent:

\begin{itemize}
  \item [(i).] $(M,\xi)$ is  non-orderable;

  \item [(ii).] There exists a contractible positive loop of
  contactomorphisms for $(M,\xi)$.

\end{itemize}

\end{prop}

This order is closely related to \textit{squeezing} properties in contact geometry \cite{EP99} as well as to the existence of bi-invariant metrics on  $\widetilde{Cont}_0(M,\xi)$ or on the space of Legendrian submanifolds \cite{CS12}.

From the beginning of the $80$'s, it is  known that 
the world of contact structures splits in two classes
with opposite behaviors. Following Eliashberg,
we say that a contact structure $\xi$ on 
$M^3$ is overtwisted if there exists an overtwisted disk
$D_{OT}\subset M$, i.e. an embedded disk which
is tangent to $\xi$ along its boundary. The overtwisted 
contact structures are flexible and classified by an adequate
h-principle \cite{EY89}. We denote $\alpha_{OT}$ a contact form for 
an overtwisted contact structure  $\xi$ defined on 
a neighborhood of an overtwisted disk. More recently, the work of Niederk\"uger \cite{KN06} and Borman-Eliashberg-Murphy \cite{BEM14} have described a similar dichotomy in the higher dimensional case. Following a suggestion of Niederk\"uger, we
say a contact structure $\xi$ is overtwisted if $(M^{2n+1},\alpha)$
contains $D^3\times D^{2n-2}(r)$ with 
 $\alpha|_{D^3\times D^{2n-2}(r)}= \alpha_{OT}-(ydx-xdy)$ for some constant $r>0$
large enough depending on the dimension of $M$ \cite{casals15}.
As in dimension three,  
Borman, Eliashberg and Murphy \cite{BEM14} have shown that
 overtwisted contact structures are  purely 
topological objects and are flexible.

On the contrary, we say $\xi$ is a tight contact structure
if it is not overtwisted. For example, the contact manifolds
$(\mathbb{S}^3,\xi_{std})$ and $(\mathbb{R}P^3,\xi_{std})$
are tight according to the fundamental result of Bennequin \cite{Be}. Similar results hold in higher dimension, where holomorphic methods give that a Liouville fillable contact structure is tight, see \cite{KN06}.

The orderability property is not shared by all contact manifolds (see the work of Albers, Frauenfelder, Fuchs and Merry \cite{AF12,AFM13,AM15} for more examples ).

\begin{thm}
\begin{itemize}
  \item [(i).] $(\mathbb{S}^3,\xi_{std})$ is non-orderable while $(\mathbb{R}P^3,\xi_{std})$ is orderable \cite{EKP06};
  \item [(ii).] There are some overtwisted contact manifolds which are non-orderable \cite{CPS14}.
\end{itemize}
\end{thm}

It is interesting to see that tight contact manifolds can be orderable or not despite their rigid nature.
At the same time we guess
overtwisted contact manifolds are non-orderable.

\begin{quest}

Are all overtwisted contact manifolds non-orderable?
\end{quest}

In order to answer the above question, we transfer the study of positive
contact isotopies to that of positive Legendrian isotopies by the  trick
of contact product. Indeed, a positive contact isotopy of $(M,\xi)$
can be lifted to a negative Legendrian isotopy of the diagonal $\Delta_{M\times M}\times \{0\}$ in the contact product
 $(M\times M\times \mathbb{R}, \alpha_1-e^s\alpha_2)$. Here $\alpha_1$ and $\alpha_2$ denote the pull-backs of $\alpha$ by the first and second projection from  $M\times M \times \mathbb{R}$ to $M$.
  The advantage is that the study of positive 
Legendrian isotopies should be  easier than
that of contact isotopies.

In that context, there is a natural question regarding positive Legendrian isotopies:
\begin{quest}\label{quest01}
Let $(M,\xi)$ be a contact manifold
and let $L_0$ and $L_1$ be Legendrian submanifolds
in $(M,\xi)$ which are Legendrian isotopic. Does there exist a positive Legendrian
isotopy connecting them?
\end{quest}

\begin{ex}\label{wave}
Let $(\mathbb{S}^2, g)$ be the $2$-sphere with 
the round metric
$g$, and let $ST^*\mathbb{S}^2$ be the space 
of contact elements on $\mathbb{S}^2$.
Denoting $S,\, N$ the poles, then the geodesic
flow of $g$ induces a positive Legendrian isotopy
$L_t$
connecting the Legendrian fibers $ST^*_N\mathbb{S}^2$
and  $ST^*_S\mathbb{S}^2$.
\end{ex}
Generally, the answer  to Question \ref{quest01}
is negative. 

\begin{thm}\label{thm0}

Let $M^n,\, n>1$ be a manifold with open 
universal cover. Then
\begin{itemize}
\item[(i).] the fibers of  $ST^*M$ are not
in a positive loop of Legendrian embeddings \cite{CFP10,CN10};
\item[(ii).] the zero-section of  $(T^*M\times \mathbb{R}, dz-ydx)$
is not in a positive loop of 
Legendrian embeddings \cite{CFP10}.

\end{itemize}

\end{thm}

We let $F$ be the front projection
 $(T^*M\times \mathbb{R},dz-ydx)\rightarrow M\times \mathbb{R}:(x,y,z)\mapsto (x,z)$. 
 For a Legendrian submanifold
 $L\subset (T^*M\times \mathbb{R},dz-ydx)$, the subset
 $L_F:=F(L) \subset M\times \mathbb{R}$ is
 the  front of 
 $L$. We usually identify $L_F$ with $L$, since
 the $y$ coordinates are given by the slopes of the front.
 In the case where $L$ and $M$ are of dimension $1$, we can replace a smooth segment of $L_F$ by a zigzag with two cusps. The zig-zag either has a $z$-shape, as in Figure~\ref{rota}, or an $s$-shape (the symmetric of Figure~\ref{rota} by the vertical axis).  
 The Legendrian submanifold obtained by this operation is 
 denoted by $S(L)$ and is called 
 a stabilization of  $L$. When we want to make it clearer, we will discrimate between the $z$-shape/positive stabilization denoted $S_+(L)$ and the $s$-shape/negative stabilization $S_-(L)$.

 We have:

\begin{prop}\cite{CFP10} \label{thmCFP}
 Let $L$ be the zero-section of  $T^*\mathbb{S}^1\times \mathbb{R}$ and
 $S(L)$ a stabilization of $L$. Then
 there  exists a loop of positive Legendrian embeddings
 based at $S(L)$.
 \end{prop}
 
For a contact manifold $(M,\xi)$ of dimension
strictly higher than three,
Murphy \cite{EM13} introduced the class of loose Legendrian submanifolds.
This is a higher dimensional generalization of    the stabilized 
 $S(L)$ in dimension three. Loose Legendrian submanifolds
 satisfy a 
h-principle discovered by  Murphy which
make them  flexible. 
The main result of this article
extends this flexible behavior.

\begin{thm}\label{thm1}
Let $(M,\xi)$ be a contact manifold of dimension $\geq 5$ and  $L\subset (M,\xi)$ be a
Legendrian submanifold. If $L$ is loose then there exists
a contractible positive loop of Legendrian embeddings based at $L$.
\end{thm}

Without the looseness assumption, F.Laudenbach has
proven that there always exist positive loops
of Legendrian immersions \cite{LF07}.

As an application of  Theorem~\ref{thm1}, we obtain a holomorphic curve free proof of the existence of tight (i.e. non overtwisted in the Borman-Eliashberg-Murphy sense \cite{BEM14}) contact structures in every dimensions. The ``hard part'' of the argument uses Theorem~\ref{thm0} whose proof relies on the existence of a generating function for a specific class of Legendrians (in that case the Legendrian fibers of the Legendrian fibration in $(\mathbb{R}^n \times \mathbb{S}^{n-1},\xi_{std} )$). 

\begin{cor}\label{tight}\cite{MNPS13}
The contact manifold $(\mathbb{R}^n \times \mathbb{S}^{n-1},\xi_{std} )$ is tight.

\end{cor}
This corollary is proved in Subsection~\ref{section:tight}. 

In the last section, we define a new partial order on certain groups $\widetilde{Cont}_0(M,\xi)$, called strong orderability, based on the transfer of an isotopy of contactomorphisms to a Legendrian isotopy of their graphs in the contact product. We then drop the graph condition to stick to Legendrian isotopies and get a (possibly) different notion than that of Eliashberg-Polterovich's \cite{EP99}.

\begin{prop}
Let $(M,\xi)$ be a contact manifold. Then
$(M,\xi)$ is strongly 
orderable if and only if 
there does not exist a 
contractible positive loop of 
Legendrian embeddings based at the  diagonal
of the contact product of $(M,\xi)$.
\end{prop}

As an example we prove that the contact manifold
$(\mathbb{S}^1,d\theta)$ 
is strongly orderable.

\bigskip

In that context, we explain the following result which was 
first suggested by Klaus Niederkr\"uger and 
also observed by
Casals and Presas.

\begin{prop}\label{overtwisted01}
Let $(M^{2n+1}, \alpha)$ be a compact overtwisted contact manifold. Then the contact
product
$(M\times M\times \mathbb{R}, \alpha_1-e^s\alpha_2)$ is also overtwisted and the diagonal $\Delta \subset (M\times M\times \mathbb{R}, \alpha_1-e^s\alpha_2)$ is loose.
\end{prop}

Therefore, according to Proposition \ref{overtwisted01}, we
have the following result:

\begin{thm}\label{thm2}
Overtwisted contact manifols are not strongly orderable.
\end{thm}

\textbf{Organisation of the paper:} In section 
$1$, we recall some basic definitions including
 Murphy's loose Legendrian embeddings. In section $2$, we give the proof
of Theorem \ref{thm1}. Finally, we prove all
the other results mentioned above in the last section.

\textbf{Acknowledgment:} This paper comes
from my thesis work. First of all,
I wish to express my infinite gratitude to
my adviser Vincent Colin for his encouragements
and countless help. I also thank Francois Laudenbach, 
Baptiste Chantraine, Claude Viterbo and Francisco
Presas sincerely for their appreciation and indispensable 
advices. I have already presented the main results
in several conferences and seminar talks, I
would like to thank the audiences for their attentions and valuable questions and suggestions.

\section{Basic definitions in Contact Geometry}

Let $L: Y\hookrightarrow(J^{1}(Y), \alpha)$
be a  smooth  Legendrian embedding. We denote 
its front map by
$L_F: Y \rightarrow Y\times \mathbb{R}$.

Given a Legendrian submanifold $L'$, there is a neighborhood $U(L')$ of $L'$ contactomorphic to a 
neighborhood of the zero section in $J^1(L',\alpha)$, according to the Weinstein neighborhood theorem. If $L$ a Legendian submanifold  close to $L'$ then  we can
talk about the front $L_F$ of $L$ in this Weinstein neighborhood.

If $\phi_t: Y \rightarrow Y\times \mathbb{R}$  is 
an homotopy of fronts (with $\phi_t(Y)$ transverse to the
$\mathbb{R}$ factor), we denote $\widetilde{\phi}_t$ its Legendrian lift and write $v_{\phi_t}$ and $v_{\widetilde{\phi}_t}$ for the corresponding generating time dependent vector fields.

\subsection{Positive Legendrian isotopies}

\begin{defn}\cite{CFP10,CN10}(\textbf{Positive Legendrian isotopy})
 Let $(M,\xi=ker\alpha)$ be a contact manifold, $L\subset M$ a Legendrian
 submanifold, $\varphi: L\times [0,1] \rightarrow M$ a Legendrian isotopy
 and let $X_t=\frac{d\varphi}{dt}$ where $t\in [0,1]$.
 We say $\varphi$ is \textbf{positive}
 if $X_t$ is transverse to $\xi$ positively, i.e.
$$\alpha(X_t)> 0.$$
Moreover,  $\varphi$ is said to be  a 
\textbf{positive loop} if in addition
$\varphi_0(L)=\varphi_1(L)$.

\end{defn}

\begin{figure}[h]

\centering\includegraphics[width=0.6\textwidth]{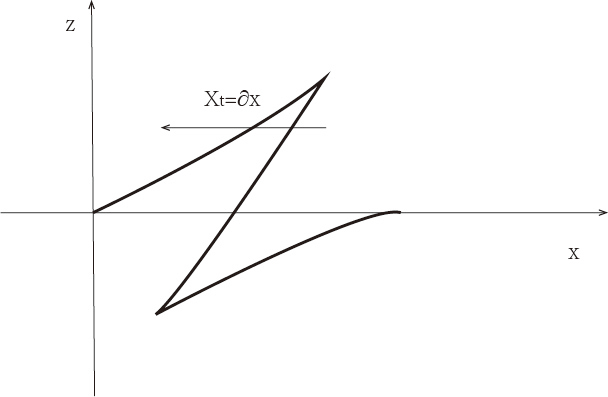}
\caption{A positive stabilized front.}\label{rota}
\end{figure}

The following remark is the starting point of our study.
\begin{rem}\cite{CFP10} \label{ps-lem1}
Let $L : \mathbb{S}^1 \hookrightarrow(J^{1}({\mathbb{S}^1}), \xi_{std})$ be
 a  Legendrian embedding whose front have positive slopes everywhere.  
 Then there exists a positive Legendrian loop 
 based at $L$.

\end{rem}

\begin{proof}
Regard $\mathbb{S}^1$ as $\mathbb{R}/\mathbb{Z}$ with coordinate $x$. Denote $Z=L_F(\mathbb{S}^1)$. On $Z$, the  slopes $\partial{z}/\partial
x>0$ are positive. See figure \ref{rota}.
Consider the vector field $X_t:= -\partial_x$
on $J^1(\mathbb{S}^1)$ and its flow $\varphi_t$. 
  
Because $\alpha(X_t)>0$ on $\varphi_t(Z)$ for every $t\in [0,1]$,  then
$\varphi_t$ is  a positive Legendrian isotopy. Since $\varphi_1=Id$, then we have a positive loop.

\end{proof}

\begin{rem}
If the front of $L$ has negative slopes everywhere, we
can choose $v=\partial_x$ so that its flow is a positive 
loop.
\end{rem}

\subsection{Loose Legendrian embeddings}

In this section, we recall Murphy's notion of loose 
Legendrian embeddings, wrinkled Legendrian 
embeddings and the idea for resolving wrinkles \cite{EM13}. For simplicity, we give the following equivalent definition of a loose Legendrian.
\begin{defn} \label{defloose2}
Let $L: Y^n\hookrightarrow (J^1(Y^n),\xi_{std})$ be a Legendrian 
embedding. Let
$\Lambda$ be a one dimensional zigzag
 and $N$ be a closed $n-1$ dimensional manifold. We say
 $L$ is loose if its front contains $\Lambda \times N$.
 In particular, it is obtained from a Legendrian $L'$ by replacing a neighborhood of $N \subset L'$ by $N$ times a zigzag. We denote $L=S_{\pm}^N (L')$, where $\pm$ stands for the $z$- or $s$-shape of the zig-zag.

\end{defn}

\begin{figure}[H]

\centering\includegraphics[width=0.7\textwidth]{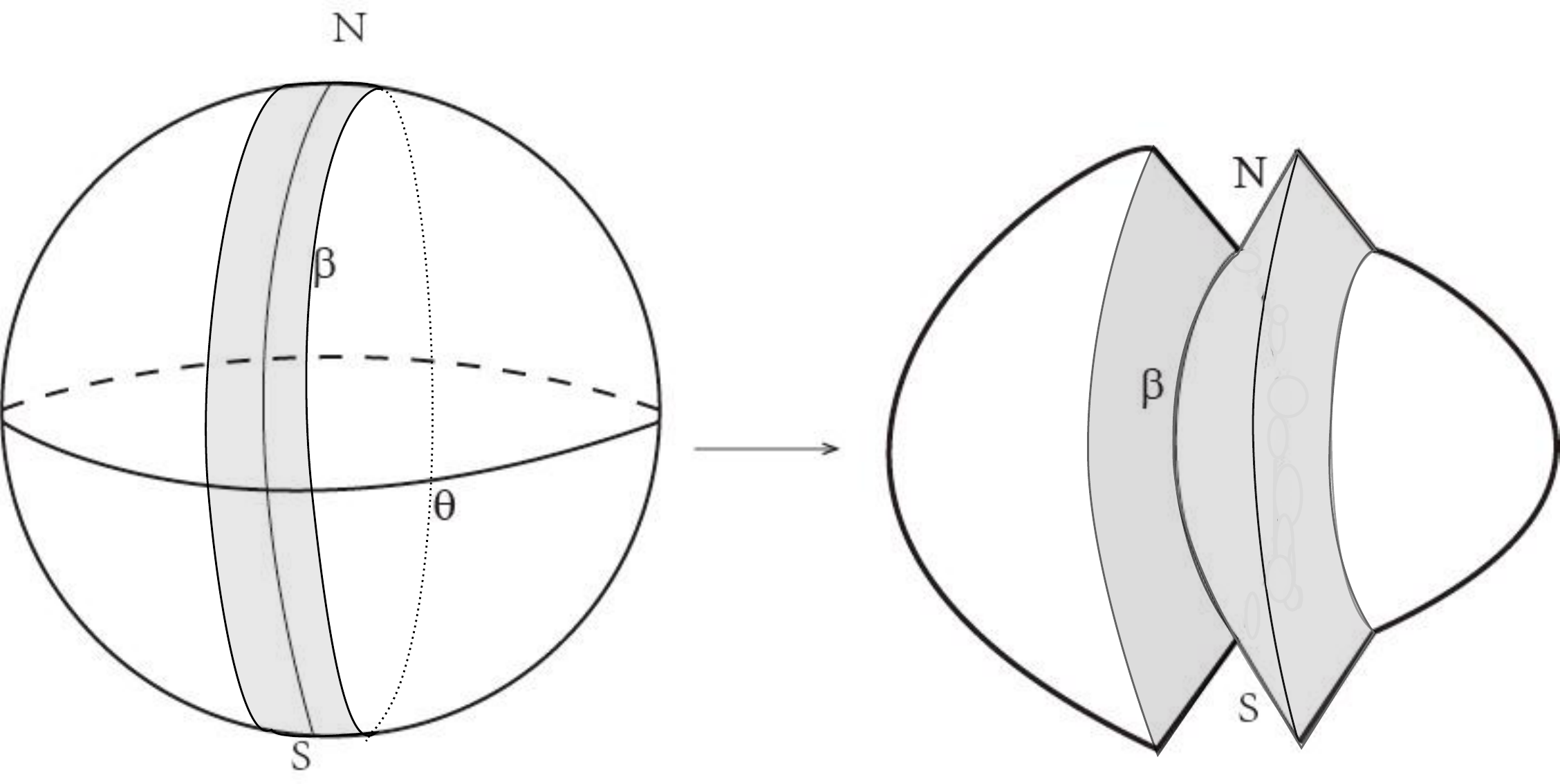}
\caption{A loose embedding of $\mathbb{S}^2$.}\label{ls}
\end{figure}

\begin{figure}[H]

\centering\includegraphics[height=2.5in]{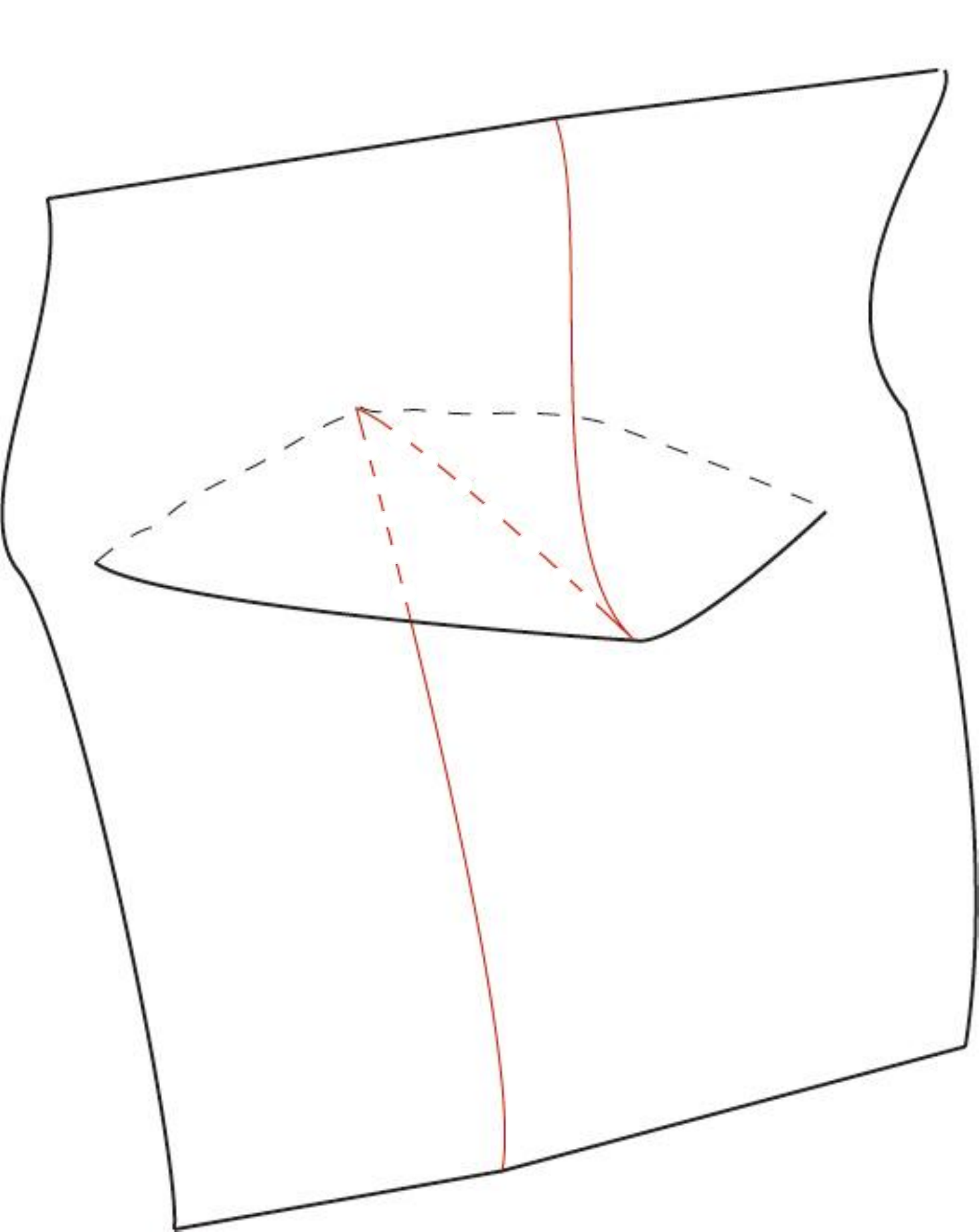}
\caption{Wrinkled embedding.}\label{Wrinkled}
\end{figure}

\begin{defn}\cite{EM11}(\textbf{Wrinkled embeddings}) 
	See figure \ref{Wrinkled}.
Let $W : \mathbb{R}^n \rightarrow \mathbb{R}^{n+1}$ be a smooth, proper map, which is a topological embedding. 
Suppose $W$ is a smooth embedding away from a finite collection
of spheres $\{S_j^{n-1} \}$. Suppose, in some coordinates near these spheres, that $W$ can be
parametrized by
\[
W(u, v) = (v, u^3-3u(1-|v|^2 ), \frac{1}{5}u^5-\frac{2}{3} u^3(1-|v|^2 )+u(1-|v|^2)^2,\]
where our domain coordinates lies in a small neighborhood of the sphere $\{|v|^2 +u^2 =
1\} \subset \mathbb{R}^n$. Then $W$ is called a 
\textit{wrinkled embedding}, and the spheres $S_j^{n-1}$ are called
the \textit{wrinkles}.
\end{defn}

 \begin{defn}\cite{EM13}(\textbf{Wrinkled Legendrians})
 Let $Y^n$ be a closed and connected manifold
  and $(M^{2n+1},\xi)$ be 
 a contact manifold.
 A wrinkled Legendrian is a smooth
 map  $L:Y \rightarrow M$, which is a topological embedding, satisfying the following properties: The image of $dL$ is contained in $\xi$ everywhere and $dL$ is full rank outside
 a subset of codimension $2$. This singular set is required to be diffeomorphic to a
 disjoint union of $(n-2)$-spheres $\{S^{n-2}_j\}$, 
 whose images are  called \textit{Legendrian wrinkles}. We assume
 the image of 
 each $S^{n-2}_j $ is contained in a Darboux chart $U_j$, so that the front projection of $L(Y)\cap U_j$
 is a wrinkled embedding, smooth outside of a compact set.
 \end{defn}

\begin{figure}[h]
\centering
\includegraphics[width=0.6\textwidth]{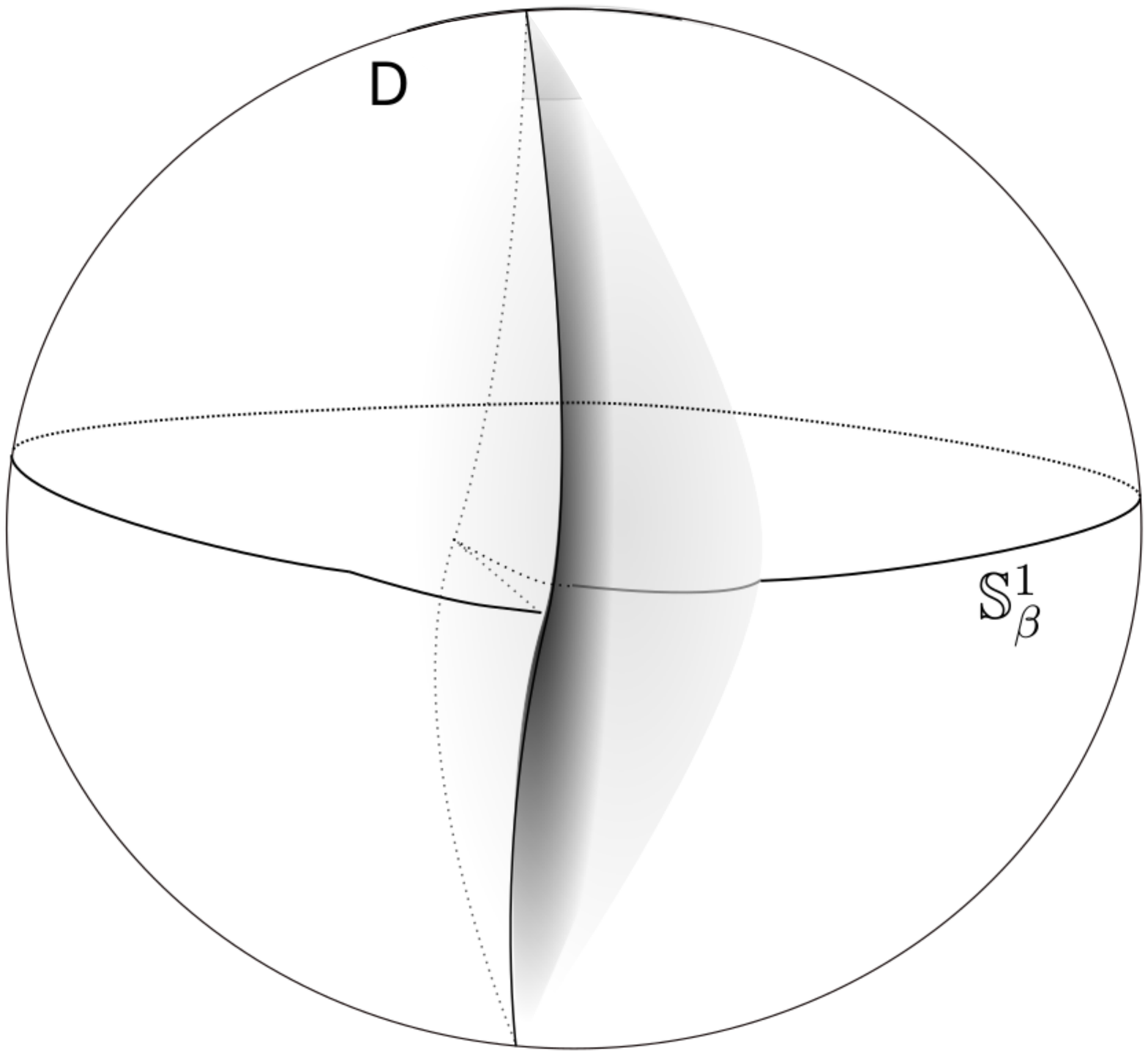}
\caption{A wrinkled sphere.}
\label{wrinkledsphere}
\end{figure}

\begin{defn}\cite{EM13}(\textbf{Twist marking})
Let $L: Y\rightarrow (M, \xi)$ be a wrinkled Legendrian 
embedding, and  $\{S^{n-2}_j\}$ be the set of 
singular spheres. Let $N \subset Y$ be a submanifold with
$\partial N =\cup_jS^{n-2}_j$.  Denote $\Phi:= L|_N$.
Then $
(\Phi, N) $ is 
called a \textit{twist marking}.
\end{defn}

\begin{rem}
We will put the $C^\infty$-topology on
the space of wrinkled Legendrian embeddings. Thus we 
can talk about a smooth family of wrinkled embeddings
$(L_t, \Phi_t, N_t)$.
\end{rem}

Given a Legendrian $L$, we
 denote $L^w$ a wrinkled Legendrian obtained by adding some wrinkles to
$L$. Given a twist marking $N$ on $L^w$ and $\eta>0$, we denote $W_{\eta,N}^{-1} (L^w)$, or $W_{\eta}^{-1} (L^w)$ when $N$ is understood, be the operation of resolving 
the wrinkles along $N$ with an $\eta$-small operation. To measure proximity, we can first perform an immersed resolution $W_0^{-1} (L^w)$, where along $N$ we incorporate a completely flat zig-zag, covering a segment times $N$. Then $W_\eta^{-1} (L^w)$ is $\eta$-$C^\infty$-close to  $W_0^{-1} (L^w)$.
This operation can be done parametrically, as summarized in the following theorem of  \cite{EM13}.

\begin{thm}\label{reswk} \cite{EM13}
Let $L^w_t$ be a smooth family of wrinkled Legendrian 
embeddings, let $(\Phi_t, N_t)$ be the twist markings.
Then there is a smooth family of  Legendrian embeddings
$L_t$, such that $L_t$ is identical to $L_t^w$ outside
of any small neighborhood of $N_t$ for all $t$.
Also, the resolution $L_t$ can be taken to be as close as we want from $L_t^w$.
\end{thm}

\begin{figure}[h]
\centering
\includegraphics[width=0.7\linewidth]{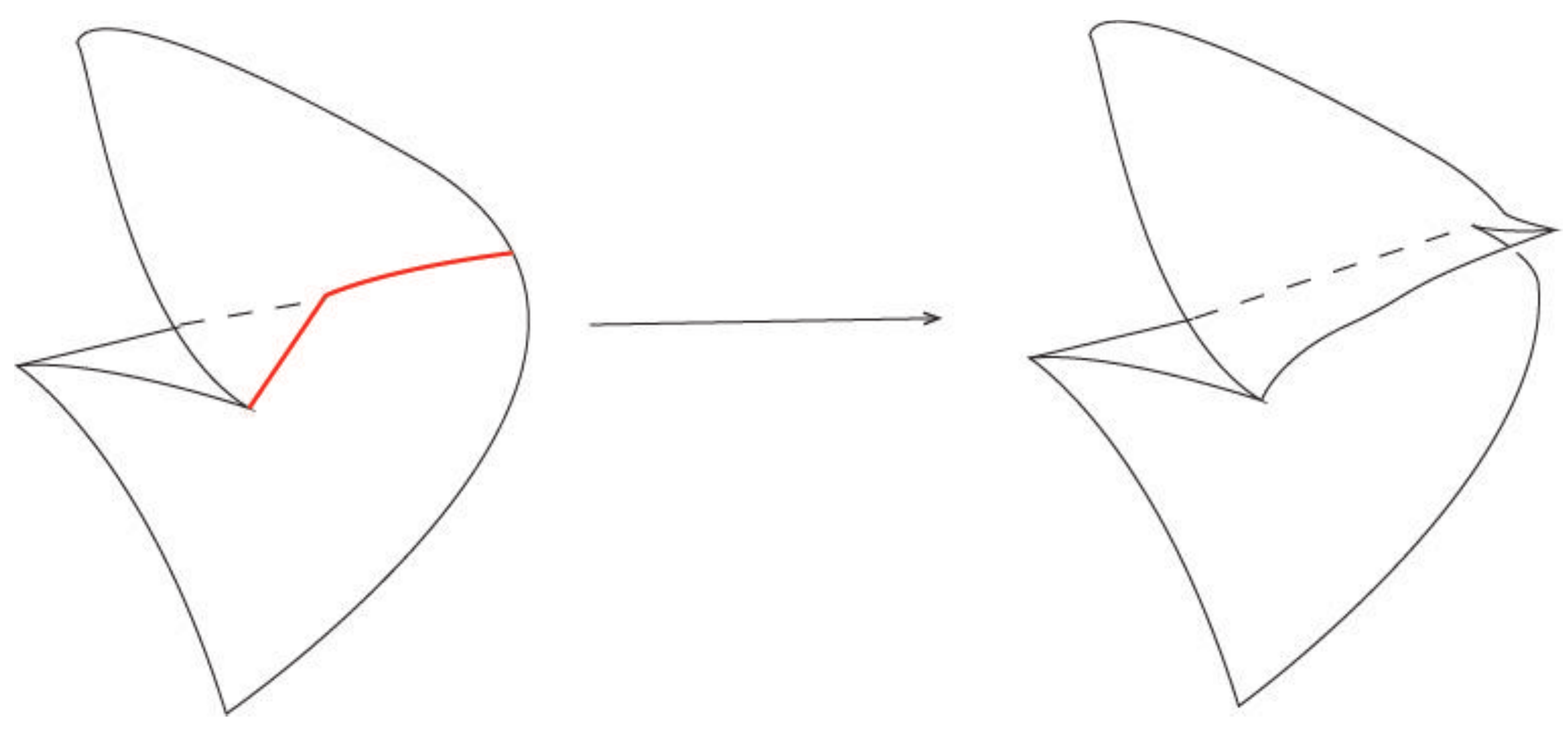}
\caption{Resolve a wrinkle.}
\label{fig:rw}
\end{figure}

\begin{thm}\label{nc}\cite{EM13}
  Let $L^n$ be a loose Legendrian and $N\subset L$ be a closed codimension $1$ submanifold of Euler characteristic $0$. If we stabilize $L$ positively and negatively along $N$, we obtain a Legendrian $S_-^N (S_+^N (L))$ which is isotopic to $L$.
\end{thm}

\begin{proof}[Sketch of proof.]
 First of all, the stabilized Legendrian  $S_-^N (S_+^N (L))$ is loose. By Murphy's $h$-principle, have to show that $S_-^N (S_+^N (L))$ is formally Legendrian isotopic to $L$. This is a consequence of the fact that $S_-^N (S_+^N (L))$ is obtained from a Murphy $N$-stabilization of $L$ by an isotopy, which doesn't change the formal Legendrian isotopy class when $\chi (N)=0$ (Proposition 2.6 in \cite{EM13}).
\end{proof}
  
In our case, we note that the stabilization operation passing from $L$ to $S_+^N (L)$ might change the formal isotopy class of the Legendrian $L$, even if $N$ has Euler characteristic zero. However, we can go back to the original formal class by stabilizing again to $S_-^N (S_+^N (L))$. This fact will be used later on in the proof of our main theorem to correct formal classes.

\section{Contractible positive Legendrian loops}
In this section we prove our main theorem \ref{thm1}
in a geometric way.
 
\begin{proof}
  We start with a loose Legendrian $L$ and work in a compact region of its standard neighborhood $J^1 (L)$.

\bigskip
  
\textbf{A. Construction of a positive loop}
 
 We first describe an elementary operation that will be applied repeatedly.
 Recall $L_F$ is the front projection of $L$. Since $L$ is the zero-section in $J^1(L)$, one can canonically identify $L_F$ with $L$.
 
 We consider a $n$-disk $D_0^n\subset L_F$, written as $D_0^{n-2}\times D^2(2)$ together with coordinates $(u,\rho, \theta)$, where $(\rho \leq 3,\theta)$ are polar coordinates on $D^2$. 
 We let $L^w$ be the wrinkled Legendrian obtained by adding 
 one wrinkled disk $D_0^w$ along the $(n-1)$-disk $D_0=\{ 1\leq \rho \leq 2 ,\theta=0\}$ to $L_F$, so that $D_0^w \subset =\{ 1\leq \rho \leq 2 \}$.
 We moreover slightly modify $L^w$ along $D_0 \times S^1 = \{ 1\leq \rho \leq 2 \}$, where $S^1$ corresponds to the $\theta$ direction, by propagating the slope of the wrinkle in the $\theta$-direction so that:
 \begin{itemize}
 \item every circle $L^w \cap \{ \rho =\rho_0 \in (1,2), u=u_0\}$, contained in the $\{(u_0,\rho_0, \theta ,z)\}$ cylinder, has a positive slope, i.e. is positively transverse to $\partial_\theta$ as in Figure~\ref{rota};
 \item $L^w$ is equal to $L_F$ away from $\{ 1\leq \rho \leq 2 \}$.
 \end{itemize}
 
 The situation is pictured in Figure \ref{wrotation}.
 \begin{figure}[h]
\centering
\includegraphics[width=0.6\linewidth]{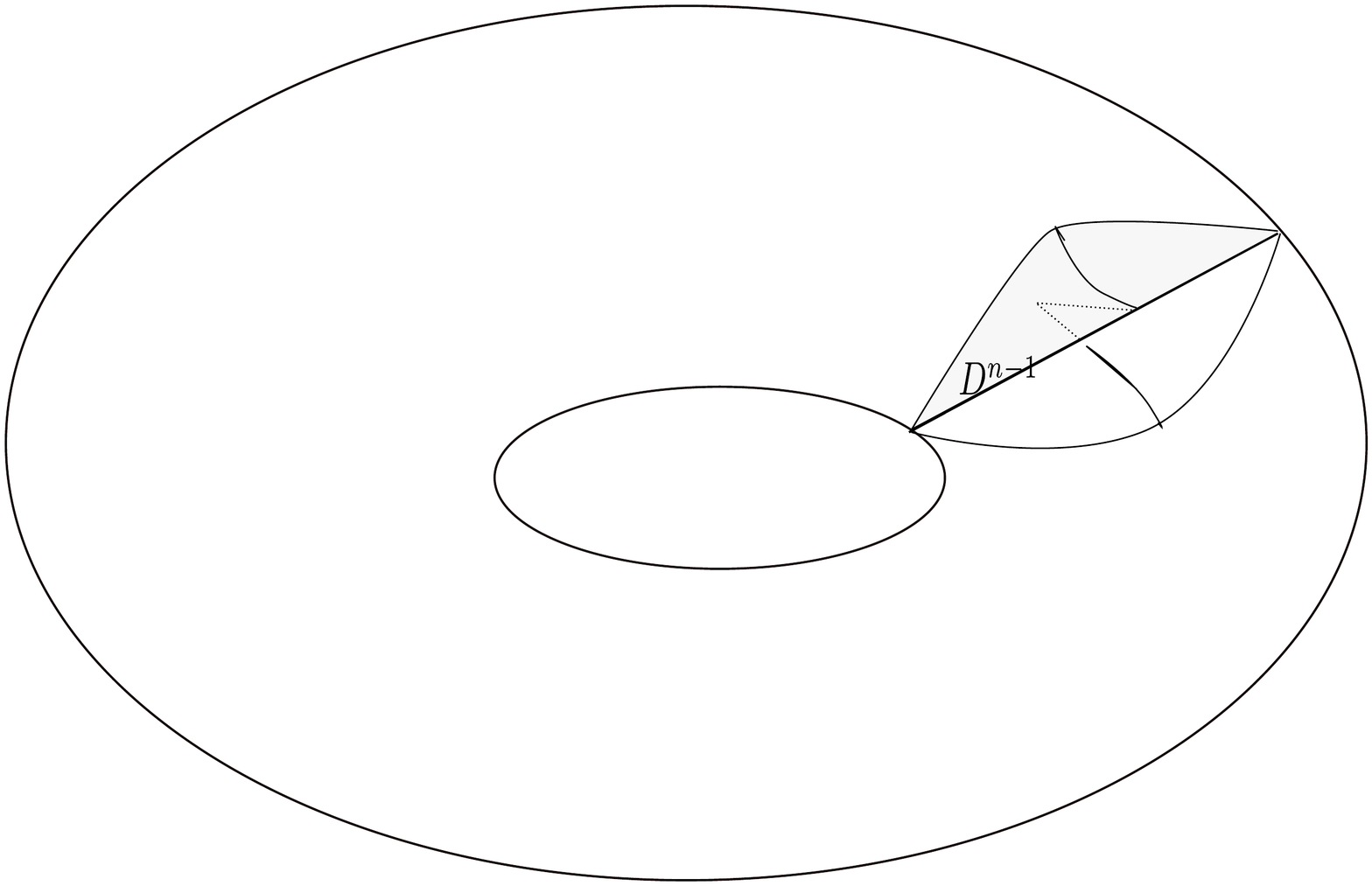}
\caption{Rotation of the wrinkle.}
\label{wrotation}
\end{figure}
We take a twist marking $N\subset D_0^n$ for $D_0^w$
so that if $N'$ is the closed submanifold $N\cup D_0 \subset L_F$ then the Euler characteristic $\chi (N')$ of $N'$ is zero.

 \smallskip

\emph{Step 1}. \textbf{We rotate the wrinkle positively in the $\theta$-direction}.

Given some constant $K_0 > 0$, we rotate positively the wrinkle in the $\theta$-direction with (large) speed $K_0 \in \mathbb{N}$: We take a 
$z$-invariant path of diffeomorphisms $\phi^0_t$ in the front such 
that $\phi^0_t(u,\rho, \theta, z)=(u, \rho, \theta-2K_0\pi t, z)$ and $\phi^0_t$ globally preserves $L^w \setminus \{ 1\leq \rho \leq 2\}$. By construction $L^w_{t,K_0} =\phi^0_t (L^w)$ is a loop of wrinkled fronts, even if $\phi^0_t$ is only a path.
Its lift 
$\widetilde{L^w_{t,K_0}}$ is a non-negative loop 
 of wrinkled Legendrians based in $L^w$. Non-negativity comes from the fact that the infinitesimal generator of the isotopy $\phi^0_t$ is either tangent (outside $\{ 1\leq \rho \leq 2\}$) or positively transverse (inside $\{ 1\leq \rho \leq 2\}$) to the front of $L^w_{t,K_0}$.
 \smallskip
 
 \emph{Step 2}. \textbf{We resolve the wrinkle}.
 We now parametrically resolve the wrinkle $\phi^0_t (D_0^w ) \subset L^w_{t,K_0}$ along the marking $\phi^0_t (N) \subset L^w_{t,K_0}$ to get a loop $W_{\eta_0}^{-1} (L^w_{t,K_0})$ of Legendrian fronts (notice that $\phi^0_1 (N) =N$). Doing so, we might introduce some negative displacement near $\phi^0_t (N)$, but, taking the size of the resolution $\eta_0$ small enough in front of $K_0$ and the slope of the circles $\{ \rho =\rho_0 \in [5/4,7/4], u=u_0\}$, we can make sure that the isotopy is still positive in the region $\{ 5/4 \leq \rho \leq 7/4 \}$.

 Here, we notice that the loop of fronts $W_{\eta_0}^{-1} (L^w_{t,K_0})$ is obtained from $L_F$ by replacing a neighborhood of
 $\phi^0_t (N')$ by the product of $N'$ with a $z$-shape. In particular, the loop of fronts $W_{\eta_0}^{-1} (L^w_{t,K_0})$ admits a parametrization by a loop of homeomorphisms (which are diffeomorphisms except at the cusps of the fronts) $\psi_{t,K_0}^0 : L \to W_{\eta_0}^{-1} (L^w_{t,K_0})$ which is constant away from $D_0^n$.
 When we lift the loop of fronts $W_{\eta_0}^{-1} (L^w_{t,K_0})$ to a loop of Legendrians $L^0_{t,K_0}$ in $J^1 (L)$, we can lift the loop of parametrizations $\psi_{t,K_0}^0$ to a loop of smooth parametrizations $\Psi_{t,K_0}^0$, which is constant away from $D_0^n$.
 The loop $L^0_{t,K_0}$ is positive along $\Psi_{t,K_0}^0 (\{ 5/4 \leq \rho \leq 7/4 \})$.
 \smallskip
 
 \emph{Step 3}. \textbf{We adjust the formal class}.
 In steps $1$ and $2$, we have constructed a loop for a stabilization $S_+^{N'} (L)$ of $L$, which might not be formally Legendrian isotopic to $L$. Following Proposition~\ref{nc}, we correct this by stabilizing again parametrically along a parallel copy of $\psi_{t,K_0}^0 (N')$ to obtain a loop of Legendrians based at $S_-^{N'}(S_+^{N'} (L))$. This second stabilization can be made small at will (in the sense that the $s$-shape is squashed) so that the positivity property of step 2 is still unchanged and we still have a loop of parametrizations, that we persist to write $\Psi_{t,K_0}^0$.
\bigskip
 
 This preparatory work been done, the proof starts from a covering of $L$ by open sets $A_i \subset D_i^n$, $i=0,\dots ,k$, of the form $S^1 \times D^{n-1} =\{ 5/4 \leq \rho \leq 7/4 \}\subset D_i^n$ as before.

 We then construct our loop by induction: by step 1,2,3, we construct a loop of Legendrians $L^0_{t,K_0} =\Psi_{t,K_0}^0 (L)$ which is positive along $\Psi_{t,K_0}^0 (A_0 )$, by rotating a wrinkle with speed $K_0$.
 We now take a loop of standard Weinstein neighborhoods $N(L^0_{t,K_0})$ of $L^0_{t,K_0}$ parametrized by $t \in S^1$, in which $L^0_{t,K_0}$ is the zero section diffeomorphic to $L$. This is given by a family of embeddings $J^1(L) \to J^1 (L)$ sending the zero section to  $L^0_{t,K_0} =\Psi_{t,K_0}^0 (L)$. These embeddings can be chosen to extend $\Psi_{t,K_0}^0$, so we still denote them $\Psi_{t,K_0}^0 : J^1(L) \to N(L^0_{t,K_0})$.

 
 We then apply steps 1,2,3 to $A_1, D_1^n$ in $J^1(L)$, though of as the source of $\Psi_{t,K_0}^0$,
 by rotating a wrinkle with (relative) speed $K_1$ and resolving it with size $\eta_1$.
 This means we are performing this step in the moving neighborhood $\Psi_{t,K_0}^0 ( J^1(L))$ of $L^0_{t,K_0}$.
 We get a loop of Legendrians $L^1_{t,K_1} =\Psi_{t,K_1}^1 (L)$ in the moving neighborhood $J^1(L)$.
Viewed in the original jet-space, we are considering the loop $\Psi_{t,K_0}^0 (L^1_{t,K_1}) = \Psi_{t,K_0}^0 (\Psi_{t,K_1}^1 (L))$.

We see  that if we take $K_1$ large enough, in particular with respect to $\eta_0$ and $K_0$, then the
loop  $\Psi_{t,K_0}^0 (\Psi_{t,K_1}^1 (L))$ becomes positive along $\Psi_{t,K_0}^0 (\Psi_{t,K_1}^1 (A_1))$ -- where it was before possibly negative.
We also have to take $\eta_1$ small enough so that the isotopy remains positive along $\Psi_{t,K_0}^0 (\Psi_{t,K_1}^1 (A_0))$ after resolving the wrinkle with size $\eta_1$.

 Precise computations are described by the following composition of speeds:

 Since
$$v_{\widetilde{\Psi}^0_{t,K_0}\circ \widetilde{\Psi}^1_{t,K_1}}(x)=v_{\widetilde{\Psi}^0_{t,K_0}}
(\widetilde{\Psi}^1_{t,K_1}(x))
+D\widetilde{\Psi}^0_{t,K_0}(v_{\widetilde{\Psi}^1_{t,K_1}
}(x)),$$
we have
$$ \alpha(v_{\widetilde{\Psi}^0_{t,K_0}\circ \widetilde{\Psi}^1_{t,K_1}}(x))
= \alpha(v_{\widetilde{\Psi}^0_{t,K_0}}(\widetilde{\Psi}^1_{t,K_1}(x)))+
(\widetilde{\Psi}^0_{t,K_0})_*\alpha(v_{\widetilde{\Psi}^1_{t,K_1}}(x)).$$

Now, we have that $\alpha(v_{\widetilde{\Psi}^0_{t,K_0}}) > -k_0$ independent of $K_1$. Moreover, since the isotopy of Legendrians is compactly supported, there exists some  $c_0>0$
independent of $K_1$
such that $\widetilde{\Psi}_{t,K_0}^{0*}\alpha =f \alpha$, where $f>c_0 >0$ in a neighborhood of the original $L$ which contains all the deformations.

We can thus see that in the neighborhood of $\Psi_{t,K_0}^0 (\Psi_{t,K_1}^1 (A_1))$ where the slope of the front is larger than some $c_1>0$, 
$\alpha(v_{\widetilde{\Psi}^1_{t,K_1}\circ \widetilde{\Psi}^0_{t,K_0}})>-k_0+c_0c_1K_1$.
Thus, for $K_1$ large enough
$\widetilde{\Psi}^0_{t,K_0}\circ \widetilde{\Psi}^1_{t,K_1}$
is positive in the neighborhood of $\Psi_{t,K_0}^0 (\Psi_{t,K_1}^1 (A_1))$.
Near $A_0$ where the loop was already positive, we do not alter positivity if the size $\eta_1$ of the resolution is small enough.

Once this is understood, it is clear that we can repeat the process until we get a loop which is positive everywhere.
At each step the rotation speed has to be higher and higher with respect to previous operations.

To conclude, we observe that we have been producing a loop  based at a loose Legendrian which is formally isotopic to $L$, and thus by Murphy's theorem~\cite{EM13} Legendrian isotopic to $L$.
  
\bigskip

\textbf{B. Contractibility}
We show that the positive loop that we have been constructing is contractible amongst Legendrian loops.

The construction was inductive on the set of annuli $(A_i)$ and thus it is enough to check that the first loop $W_{\eta_0}^{-1} (L^w_{t,K_0})$
is homotopic to a constant loop. 

We first treat the case when the dimension of the Legendrian $L$ is greater of equal to $3$.
Define
$\phi^0_{s,t}$ such that $\phi^0_{s,t}(u,\rho, \theta, z)=(u, \rho-s, \theta-2K_0\pi t, z)$. We can see that $\phi^0_{s,t} (L^w)$ is a homotopy from $\phi^0_{t,K_0}(L^w) $ to 
$\phi^0_{1,t}(L^w)$ which is a loop of rotation of a wrinkled disk $D^w$ around some point, says $x_0$.
Up to homotopy, the wrinkled disk $\phi^0_{1,t}(D^w)$ is completely determined by its normal vector in $L$ at $x_0$, and thus by a map $S^1 \to S^{n-1}$. Since $n\geq 3$, this map is homotopic to a point and thus we can deform our loop of wrinkled Legendrians to a constant loop. Moreover, this homotopy can be extended to a homotopy of twist markings from the original loop of twists markings to a constant loop. Resolving parametrically the markings, we get a homotopy from $W_{\eta_0}^{-1} (L^w_{t,K_0})$ to a constant loop.
The extra stabilization of step 3 to fix the formal isotopy class enters the same scheme and can be also homotoped to a constant operation.
This concludes the proof.

The case when the dimension of $L$ is two follows the same scheme, except that we homotope the loop of resolved wrinkles $\Psi^0_{t,K_0} (U)$, where $U$ is a circle times a $z$-shape segment, to a constant annuli around the circle $\{ \rho=1\}$.
 \end{proof}

\section{Applications}
In this chapter, we give some applications of 
our main theorem. First, we reprove tightness
of $(\mathbb{S}^{n-1}\times \mathbb{R}^n,\xi_{std})$. Second, we define
a partial order on the universal cover
$\widetilde{Cont}_0(M,\xi)$ of the
identity component of the group of contactomorphisms of a contact manifold 
$(M,\xi)$ and prove that overtwisted contact structures are not orderable.

\subsection{Tightness of $(\mathbb{S}^{n-1}\times \mathbb{R}^n,\xi_{std})$} \label{section:tight}

In this section we prove Corollary \ref{tight}.
A similar proof for $\mathbb{S}^1\times \mathbb{R}^2$ was given in \cite{CFP10}.

\begin{proof}
Assume  $(\mathbb{S}^{n-1}\times \mathbb{R}^n,\xi_{std})$ is overtwisted,
and $D_{OT}\subset  (\mathbb{S}^{n-1}\times \mathbb{R}^n,\xi_{std})$ is an overtwisted disk.
Denote $\pi:\mathbb{S}^{n-1}\times \mathbb{R}^n
\rightarrow \mathbb{R}^n$ the projection.
There exists some point $x\in \mathbb{R}^n$ such
that the fiber $\pi^{-1}(x)\cap D_{OT}=\emptyset$.
According to \cite{casals15}, the fiber 
$\pi^{-1}(x)$ is loose.
Thus, there exists a positive loop based at it
by Theorem \ref{thm1}. That contradicts Theorem
\ref{thm0}. Therefore, the manifold $(\mathbb{S}^{n-1}\times \mathbb{R}^n,\xi_{std})$ is tight. 
\end{proof}

\subsection{Positive loops and orderings}

\begin{defn}
Given a contact manifold $(M,\alpha)$, the
manifold
 $(\Gamma_M, \tilde{\alpha})=(M\times
M\times \mathbb{R},\alpha_1 -e^s\alpha_2)$ is called
a contact product. Here
$\alpha_i=\pi_i^*\alpha$ where $\pi_i$ project $\Gamma_M$ to the $i$-th
factor.
The Legendrian submanifold of $(\Gamma_M, \tilde{\alpha})$
 $\Delta=\{(x, x, 0)\}$ is called the diagonal.
\end{defn}
 
The contact product $\Gamma_M$ is a special
case of a contact fibration. We
recall the definition from \cite{presas07}.

\begin{defn}
Let $(E,\xi=ker\alpha)$ be a contact manifold,
and $E\longrightarrow B$ is a fibration with fiber
$F$. Then  $(E,\xi=ker\alpha)\longrightarrow B$
is called a contact fibration if  $(F,\alpha|_F)$ is a
contact manifold.
Let $(E,\xi=ker\alpha)\longrightarrow B$ be contact fibration
we say that the horizontal distribution $H=(TF\cap\xi)^{\perp d\alpha}$ 
is the contact connection associated to the fibration.
\end{defn}

\begin{rem}
The horizontal distribution depends on the contact form $\alpha$.
\end{rem}

The connection defined above has the following properties:
\begin{prop}\cite{presas07}
For a path $\gamma :[0,1]\rightarrow B$, the monodromy $m_\gamma : F(\gamma(0))\rightarrow F(\gamma(1))$
 induced by $\gamma$ is a contactomorphism.

\end{prop}

\begin{cor}\label{hl}
Let $\phi\in Diff_0(B)$. Then it lifts to 
a contactomorphism $\widetilde{\phi}$.
\end{cor}

Note that $\Gamma_M$ is a contact fibration with $F=(M,\alpha)$ and $B=M\times \mathbb{R}$.

We now explain the following result which was 
first suggested by Klaus Niederkr\"uger and 
also observed by
Casals and Presas.
\begin{prop}\label{overtwisted}
Let $(M^{2n+1}, \alpha)$ be a compact overtwisted contact manifold and let $(\Gamma_M,\tilde{\alpha})$
be the associated contact product. Then
$(\Gamma_M,\tilde{\alpha})$ is also overtwisted and the diagonal $\Delta \subset \Gamma_M$ is loose.
\end{prop}

\begin{proof}

We  apply the overtwisted criterion from \cite{casals15}.
If $\lambda=ydx-xdy$,
it is  enough to construct a higher dimensional
overtwisted ball
$D=(B_{OT}^{2n+1}\times D^{2n+2}(r),  \alpha_{OT}-\lambda)\subset (\Gamma_M,\tilde{\alpha})$ for some $r$
large enough,
 such that $D$ does not intersect $\Delta$.
 
 Let $\mathbb{S}^{2n+1}=\{(x,y)\mid 
 x^2+y^2=1\}$ with its standard
 contact form $\alpha_{std}$, and 
 let $\varphi_0: \mathbb{S}^{2n+1} \times 
 \mathbb{R} \rightarrow \mathbb{R}^{2n+2},
  (x,y,s)\mapsto (e^sx,e^sy)$. Note that
  $\varphi_0^*\lambda=\alpha_{std}$.
We take a Darboux ball $B\subset (M,\alpha)$
and we regard it as a subset of $(\mathbb{S}^{2n+1},\alpha_{std})$.
Then we can construct a  contact embedding 
 $\varphi : (M\times B\times \mathbb{R},\tilde{\alpha})\hookrightarrow
 (M\times \mathbb{R}^{2n+2},\alpha_1-\lambda)$
by the following series of contact embeddings 
$$(M\times B\times \mathbb{R},\tilde{\alpha})
\stackrel{i}{\hookrightarrow} (M\times \mathbb{S}^{2n+1} \times \mathbb{R},
\alpha_1-e^s\alpha_{std})\stackrel{id \times \varphi_0}{\longrightarrow} (M\times \mathbb{R}^{2n+2},\alpha_1-\lambda).$$
Let $B_{OT}^{2n+1}\subset M $ be a 
overtwisted ball, then $D_0=(B_{OT}^{2n+1}\times \mathbb{D}^{2n+2}(r),\alpha_1-\lambda)$ is the overwisted ball in $(M \times \mathbb{R}^{2n+2},\alpha_1-\lambda)$.
We can move $D_0$ away from $\varphi(\Delta)$ by 
 Corollary \ref{hl}. More precisely, we 
  take the vector
 field $V=2r\partial_x+ 2r\partial_y$ on
 $\mathbb{R}^{2n+2}$, then lift it  to a contact 
 vector field 
 $V'=V+2r(y-x)R_\alpha$ on $M\times \mathbb{R}^{2n+2}$
 where $R_\alpha$ is the Reeb vector field of 
 $(M,\alpha)$. Let $\phi_t$ be the contact 
 isotopy of $V'$. 
 Denote $C=\{rx\mid x\in B,r>0\}$  the cone defined
 by $B$.
  Then $D_1=\phi_1(D_0)\subset M\times (C\setminus \{0\})=\varphi(M\times B)$ 
  does not
  intersect $\varphi(\Delta)$. Therefore
  $D=\varphi^{-1}(D_1)$ is one of the
  overtwisted balls we want.
  
\end{proof}

\begin{cor}\label{MC01}
Let $(M, \alpha)$ be a compact overtwisted contact manifold and $(\Gamma_M,\tilde{\alpha})$
the contact product. Then there exists a 
positive loop of Legendrian embeddings based at 
$\Delta$.
\end{cor}

Let $Leg(M, \Gamma_M)$ be the set of 
Legendrian embeddings $M\hookrightarrow (\Gamma_M,\widetilde{\alpha})$.
Given $\phi\in Cont_0(M,\xi=ker\alpha)$ 
with $\phi^*\alpha=e^{g(x)}\alpha$, 
it induces a contactomorphism  
$$\bar{\phi}(x,y,s):=\left(x, \phi(x),s-g(y)\right)$$
on $(\Gamma_M,\widetilde{\alpha})$. We denote 
$gr(\phi)=\bar{\phi}|_\Delta$
which is  in  $ Leg(M,\Gamma_M)$.
In fact, given a positive contact isotopy $\phi_t$,
we can see that  $gr(\phi_t)$ is
a negative Legendrian isotopy. Therefore, we would 
like to transfer the study of positive contact
isotopies to that of negative Legendrian isotopies.

\begin{defn}
Let $f=[f_t] $ and $g=[g_t]$  be  two elements in 
$\widetilde{Cont}_0(M,\xi)$. We say
$f\succeq g$ if there exists a  non-positive 
path $L_t\in Leg(M,\Gamma_M)$
from $gr(g_1)$ to 
 $gr(f_1)$ and $gr(g_t)*L_t$ is  homotopic to $gr(f_t)$. The space $\widetilde{Cont}_0(M,\xi)$
and $(M,\xi)$ are said to be \textbf{strongly orderable} if $\succeq$ defines a 
partial order \footnote{in the sense of a partial order on sets} on it. Otherwise,
they are said to be non strongly orderable.
\end{defn}

\begin{rem}
Let $C$ be the set generated by all the homotopy classes
of non-positive paths in $Leg(M,\Gamma_M)$. Then 
$f\succeq g$ equals to $gr(g^{-1}f)\in C$. Given
 $[L_t]\in C$ and  $\phi\in Cont_0(M,\xi)$, then
 we have $[\bar{\phi}L_t]\in C$. Therefore,
the order $\succeq$ is left invariant,
that is to say, given $f$ and $g$ in
 $\widetilde{Cont}_0(M,\xi)$, if $f\succeq g$, then 
 $hf\succeq hg$ for all $h\in \widetilde{Cont}_0(M,\xi) $. Because if $L_t$ is a non-positive
 path  from $g_1$ to $f_1$, then $\bar{h}_1L_t$ is 
 a non-positive path from $h_1g_1$ to $h_1f_1$.
\end{rem}

\begin{prop}
Let $(M,\xi)$ be a contact manifold. Then
$(M,\xi)$ is strongly orderable if and only if
there does not exist a contractible  negative 
loop of Legendrian embeddings based at $\Delta$.
\end{prop}

\begin{proof}
Let $f=[f_t],\, g=[g_t]$ and $h=[h_t]$  be elements in $\widetilde{Cont}_0(M,\xi)$. The relation $\succeq$
is reflective, since
we have $f\succeq f$ by
the definition of $\succeq$. If there are two non-positive
 paths $L_t^1$ from $gr(g_1)$ to $gr(f_1)$ and $L_t^2$ from
$gr(h_1)$ to $gr(g_1)$, then $L^2_t*L^1_t$ is a non-positive
path from $gr(h_1)$ to $gr(f_1)$. Thus, the relation is 
transitive.
Now we check the antisymmetry of $\succeq$.
According to \cite{CN13}[Propostion 4.5], the existence 
of contractible non-positive non-trivial loop of  Legendrian embeddings 
is equivalent to the existence of contractible negative 
loop of  Legendrian embeddings.
Thus,  for any $f\neq 1$, on one hand,
 if there does not exist any negative loop based at $\Delta$, we can
not find a non-negative path $L^1_t$ and
a non-positive path $L^1_t$ in the homotopy
class of $gr(f_t)$ at the same times. Otherwise,
$L^1_t*L^2_t$ would be a contractible non-negative loop.
On the other hand, if there exists a non-positive 
loop $f_t$ based at $\Delta$, then $f_{1/2}\succeq 1$ and $1\succeq f_{1/2}$. That means $(M,\xi)$ is 
not strongly orderable.
\end{proof}

Our definition is stronger than that of \cite{EP99},
since we do not require the path of Legendrian 
embeddings $\widetilde{\phi}_t$
 to be graphical for all $t$.

\begin{cor}
Let $(M,\xi)$ be a contact manifold. If $(M,\xi)$ is strongly orderable, then
it is orderable.
\end{cor}

A contact manifold which is not strongly orderable is said to be weakly non-orderable.
Immediately, according to Proposition 
\ref{overtwisted} and Corollary \ref{MC01},
 we deduce theorem \ref{thm2} saying that 
overtwisted contact manifols are weakly non-orderable.

We have the  following example of strong
orderability.

\begin{thm}
$(\mathbb{S}^1,\xi_{std})$ is strongly orderable.
\end{thm}
\begin{proof}
Denote $d\theta$ the standard contact form
for $\mathbb{S}^1$. We have a contactomorphism
$\varphi: (\Gamma_{\mathbb{S}^1},d\theta_1-e^sd\theta_2)
\rightarrow (\mathbb{S}^1\times T^*\mathbb{S}^1,dz-ydx), (\theta_1,\theta_2,s)\mapsto (z=\theta_1-\theta_2,
x=\theta_2,y=e^s-1)$ such that
$\varphi(\Delta)$ is the zero-section.
Assume there exists a contractible positive loop
based at the zero-section of $(\mathbb{S}^1\times T^*\mathbb{S}^1,dz-ydx)$, 
then it lifts to a positive loop based at the zero-section of $(\mathbb{R}^1\times T^*\mathbb{S}^1,dz-ydx)$. However, such loops do
not exist according to \cite{CFP10} (notice this is not a trivial result).
Thus $(\mathbb{S}^1,\xi_{std})$ is strongly orderable.
\end{proof}

\begin{quest}
Is $(\mathbb{R}P^3,\xi_{std})$ strongly orderable?
\end{quest}

\bibliographystyle{amsalpha}
\bibliography{paper}

\end{document}